\documentclass[preprint,1p]{elsarticle}

\makeatletter
 \def\ps@pprintTitle{%
 	\let\@oddhead\@empty
 	\let\@evenhead\@empty
 	\def\@oddfoot{\footnotesize\itshape
 		{} \hfill\today}%
 	\let\@evenfoot\@oddfoot
 }
\makeatother
\usepackage[unicode]{hyperref}

\usepackage{latexsym}
\usepackage{indentfirst}
\usepackage{amsxtra}
\usepackage{amssymb}
\usepackage{amsthm}
\usepackage{amsmath}
\usepackage{mathrsfs} 

\usepackage{xcolor}
\usepackage{color}

\usepackage{amsfonts}

\usepackage[capitalise]{cleveref}

\newtheorem{theor}{Theorem}[section]
\newtheorem{prop}[theor]{Proposition}
\newtheorem{cor}[theor]{Corollary}
\newtheorem{lemma}[theor]{Lemma}
\theoremstyle{definition} 
\newtheorem{defin}[theor]{Definition}
\newtheorem{rem}[theor]{Remark}

\newtheorem{ex}{Example}
\newtheorem*{conv*}{Convention}
\newtheorem{que}{Question}

\DeclareMathOperator{\Sym}{Sym}
\DeclareMathOperator{\id}{id}

\DeclareMathOperator{\Ret}{Ret}
\DeclareMathOperator{\mpl}{mpl}
\DeclareMathOperator{\Fix}{Fix}
 \DeclareMathOperator{\Dis}{Dis}

\newcommand{\seq}[1]{\left<#1\right>}

\begin{document}

\begin{frontmatter}
	\title{Simplicity and finite primitive level of indecomposable set-theoretic  solutions of the Yang-Baxter equation
\tnoteref{mytitlenote}
}
	\tnotetext[mytitlenote]{This work was partially supported by the Dipartimento di Matematica e Fisica ``Ennio De Giorgi'' - Università del Salento. The second author was partially supported by the ACROSS project ARS01\_00702. The authors are members of GNSAGA (INdAM).}
	\author[unile]{M.~CASTELLI}
	\ead{marco.castelli@unisalento.it}
		\author[unile]{M.~MAZZOTTA}
	\ead{marzia.mazzotta@unisalento.it}
	\author[unile]{P.~STEFANELLI}
	\ead{paola.stefanelli@unisalento.it}
	\address[unile]{Dipartimento di Matematica e Fisica ``Ennio De Giorgi''
		\\
		Universit\`{a} del Salento\\
		Via Provinciale Lecce-Arnesano \\
		73100 Lecce (Italy)}

\begin{abstract}
This paper aims to deepen the theory of bijective non-degenerate set-theoretic solutions of the Yang-Baxter equation, not necessarily involutive, by means of q-cycle sets. We entirely focus on the finite indecomposable ones among which we especially study two classes of current interest: the simple solutions and those having finite primitive level. 
In particular, we provide two group-theoretic characterizations of these solutions, involving their permutation groups. Finally, we deal with some open questions.
\end{abstract}

\begin{keyword}
\texttt{set-theoretic solution\sep Yang-Baxter equation\sep q-cycle set\sep cycle set\sep dynamical extension}
\MSC[2020] 
16T25\sep 81R50 \sep 20N02 \sep 20E22
\end{keyword}
\end{frontmatter}

\section*{Introduction}
A \emph{set-theoretic solution of the Yang-Baxter equation} on a non-empty set $X$ is a pair $\left(X,r\right)$, where 
$r:X\times X\to X\times X$ is a map such that the relation
\begin{align*}
\left(r\times\id_X\right)
\left(\id_X\times r\right)
\left(r\times\id_X\right)
= 
\left(\id_X\times r\right)
\left(r\times\id_X\right)
\left(\id_X\times r\right)
\end{align*}
is satisfied.  
The paper by Drinfel'd \cite{drinfeld1992some}, a milestone in quantum group theory, moved the interest of several researchers for finding solutions of this equation in the last thirty years.
Writing a solution $(X,r)$ as $r\left(x,y\right) = \left(\lambda_x\left(y\right)\rho_y\left(x\right)\right)$, with
$\lambda_x, \rho_x$ maps from $X$ into itself, for every $x\in X$, we say that $(X, r)$ is \emph{left non-degenerate} if $\lambda_x\in \Sym_X$, \emph{right non-degenerate} if $\rho_x\in \Sym_X$, for every $x\in X$,  \emph{non-degenerate} if it is both left and right non-degenerate. Moreover, $(X,r)$ is \emph{involutive} if $r^2=\id_{X\times X}$. 

Over the years, the involutive non-degenerate solutions have been widely studied starting from the seminal papers by Gateva-Ivanova and Van den Bergh \cite{gateva1998semigroups}, Gateva-Ivanova and Majid \cite{gateva2008matched}, and Etingov, Schedler, and Soloviev \cite{etingof1998set}. At the same time, Lu, Yan, and Zhu \cite{LuYZ00} and Soloviev \cite{So00} dealt with bijective non-degenerate solutions, not necessarily involutive.

A useful strategy for determining all the involutive solutions arises in \cite{etingof1998set} and consists in finding those that can not be deconstructed into other ones, the so-called indecomposable solutions. In this context, several authors introduced useful tools for determining and classifying this type of solutions (see, for example, \cite{cacsp2018,capiru2020,JePiZa20x,rump2020one,Ru220,Ru20,smock}). Later, in  \cite{CeSmVe19,EtSoGu01}, the notion of indecomposability was extended to solutions that are not necessarily involutive. Technically, a bijective non-degenerate  solution $\left(X,r\right)$ is said to be \emph{decomposable} if there exists a partition $\{X_1,X_2\}$ of $X$ such that $r\left(X_i \times X_j\right) = X_j\times X_i$, for all \mbox{$i, j\in \{1,2\}$}; otherwise $(X,r)$ is  called \emph{indecomposable}.
Etingof, Soloviev, and Guralnick \cite[Lemma 2.1]{EtSoGu01} exclusively characterized indecomposable solutions in group-theoretic terms. In this paper, we make explicit that a solution $(X,r)$ is indecomposable if and only if its \emph{permutation group}, i.e.,  the group
$$
\mathcal{G}\left(X,r\right)=\seq{\lambda_x, \eta_x \, |\, x\in X},
$$ 
acts transitively on $X$, with $\eta_x\left(y\right):=\rho_{\lambda^{-1}_y\left(x\right)}\left(y\right)$, for all $x,y\in X$.
In general, although the significant results obtained until now, finding and classifying all the indecomposable solutions is rather difficult.

Among involutive solutions, Vendramin first approached to the class of finite simple ones in \cite{vendramin2016extensions}. More recently, Ced{\'o} and Okni{\'n}ski  \cite{CeOk20x} have been introduced an equivalent definition in the indecomposable finite case; namely, an involutive solution $(X,r)$ is said to be \emph{simple} if $|X|>1$ and, for every epimorphism of solutions $f:(X, r)\to (Y,s)$, either $f$ is an isomorphism or $|Y|=1$. 
In \cite[Sections 4 and 5]{CeOk20x}, several examples of involutive simple solutions can be found.
This study was mainly motivated by the fact that the unique finite simple solutions known until then were two instances of order $4$ (see \cite[Example 2.11]{vendramin2016extensions} and \cite[Example 9]{cacsp2018}) and the \emph{primitive} ones, i.e., solutions $(X, r)$ for which $\mathcal{G}(X,r)$ acts on $X$ as a primitive group. 
In this context, motivated by a question posed by Ballester-Bolinches in \cite{oberw}, Ced{\'o}, Jespers, and Okni{\'n}ski have been recently classified the finite primitive involutive solutions in \cite[Theorem 3.1]{CeJeOk20x} and, in particular, they showed that they are all of prime order.
\medskip

The main aim of this paper is to study finite indecomposable bijective non-degenerate solutions, not necessarily involutive. 
\smallskip

\noindent \textbf{Convention.} Hereinafter, we briefly call \emph{solution} any finite bijective non-degenerate set-theoretic solution of the Yang-Baxter equation.
\smallskip

To our purposes, we fully exploit the existing one-to-one correspondence between solutions and regular q-cycle sets, algebraic structures introduced by Rump \cite{rump2019covering}. 
Specifically, a non-empty set $X$ endowed with two binary operations $\cdot$ and $:$ is said to be a \emph{q-cycle set} if the map $\sigma_x:X\rightarrow X, y\mapsto x\cdot y$ is bijective, for every $x\in X$, and the following conditions
\begin{align}
\left(x\cdot y\right)\cdot \left(x\cdot z\right) 
&= \left(y:x\right)\cdot \left(y\cdot z\right)\label{ug1}\tag{q1}\\
\left(x:y\right):\left(x:z\right) 
&= \left(y\cdot x\right):\left(y:z\right)\label{ug2}\tag{q2}\\
\left(x\cdot y\right):\left(x\cdot z\right) 
&= \left(y:x\right)\cdot \left(y: z\right)\label{ug3}\tag{q3}
\end{align}
hold, for all $x,y,z\in X$. Besides, $X$ is \emph{regular} if the map $\delta_x:X\rightarrow X, y\mapsto x:y$ is bijective, for every $x\in X$; \emph{non-degenerate} if $X$ is regular and the squaring maps, i.e., the maps $\mathfrak{q}$ and $\mathfrak{q'}$ from $X$ into itself given by
$\mathfrak{q}\left(x\right):=x\cdot x$ and $\mathfrak{q'}\left(x\right):=x: x$,
for every $x \in X$, are bijective. In \cite[Theorem 5]{CaCaSt20}, it is proved that every finite regular q-cycle set $X$ is non-degenerate.
Thus, if $X$ is a finite and regular q-cycle set, then $\left(X,r\right)$ is a solution, where $r:X\times X\to X\times X$ is the map given by 
\begin{align*}
r\left(x,y\right)= \left(\sigma_{x}^{-1}\left(y\right),\delta_{\sigma_{x}^{-1}\left(y\right)}\left(x\right)\right),
\end{align*}
for all $x,y\in X$. 
Vice versa, if $\left(X,r\right)$ is a solution, set
\begin{align*}
  x\cdot y:=\lambda_x^{-1}\left(y\right) 
  \qquad
  \text{and}
  \qquad
  x:y:=\rho_{\lambda_y^{-1}\left(x\right)}\left(y\right), 
\end{align*}
for all $x,y\in X$, then $X$ is a regular q-cycle set
(cf. \cite[Proposition 1]{rump2019covering}).
Evidently, if $X$ is a q-cycle set such that $\cdot$ and $:$ coincide, then $X$ is a \emph{cycle set}. Cycle sets were introduced by Rump in \cite{rump2005decomposition} and rather investigated (see, for instance, \cite{ cacsp2018, cacsp2019,CaCaSt20x,capiru2020,rump2016quasi,vendramin2016extensions}) for their one-to-one correspondence with left non-degenerate involutive solutions.

To study indecomposable solutions we look at indecomposable regular q-cycle sets. In light of our convention, throughout our treatise, it is clear that every q-cycle set will be finite even if it will not be specified. As one can expect, a regular q-cycle set $X$ is indecomposable if and only if its permutation group
$\mathcal{G}\left(X\right)=\seq{\sigma_x, \delta_x \, |\, x\in X}$
acts transitively on $X$.\\ 
As a first result, we show that any indecomposable q-cycle set $X$ with $|X| > 1$ and having regular permutation group $\mathcal{G}(X)$ is retractable.  In addition, if $\mathcal{G}(X)$ is also abelian, then $X$ is multipermutational.  In this regard, given a regular q-cycle set $X$, it is \emph{retractable} if $|X|=1$ or there exist two distinct elements $x, y\in
X$ such that $\sigma_x=\sigma_y$ and $\delta_x=\delta_y$. 
Besides, it is possible to consider a congruence on $X$, the \emph{retract relation} $\sim$ (see
\cite[Definition 1]{CaCaSt20}), which ensures that the quotient $\Ret\left(X\right):=X/\sim$ is a q-cycle set. Moreover, for a \emph{multipermutational} q-cycle set $X$ of level $n$, we mean that $n$ is the minimal non-negative integer such that $|\Ret^n\left(X\right)|=1$, where $\Ret^0\left(X\right):=X$ and $\Ret^i\left(X\right):=\Ret\left(\Ret^{i-1}\left(X\right)\right)$, for every $i>0$.

The core of this work is a description of indecomposable q-cycle sets in terms of dynamical extensions, a method to construct new families of q-cycle sets, already developed in \cite{CaCaSt20}. Using this tool, we construct examples of indecomposable solutions that are not involutive. As another application, we prove that the permutation group $\mathcal{G}(X)$ associated to any indecomposable retractable q-cycle set $X$ always acts imprimitively  on $X$, whenever $X$ has not prime size. Furthermore, motivated by a recent paper of Gateva-Ivanova \cite{Ga21}, we focus our attention on indecomposable square-free q-cycle sets. We recall that a non-degenerate q-cycle set is \emph{square-free} if $\mathfrak{q} = \mathfrak{q'} = \id_X$. At first, we  provide  a  structure
theorem of indecomposable square-free q-cycle sets. Moreover, referring to \cite[Question 9.6(7)]{Ga21}, we exhibit a family of these q-cycle sets that are not self-distributive.

Another aspect we deal with is that of the simplicity of q-cycle sets, whose definition is consistent with that given for involutive solutions. Initially, we show that any finite simple q-cycle set $X$ with $|X| > 1$ is indecomposable. 
Furthermore, we describe in group-theoretic terms all the finite simple regular q-cycle set $X$, involving a special group, namely, the \emph{displacement group} of $X$, given by
\begin{align*}
   \Dis(X)= \seq{\sigma_x^{-1}\sigma_y,\delta_x^{-1}\delta_y\ |\ x,y\in X}.
\end{align*}
Note that, in the finite case, our results on the displacement group include those given in \cite{bon2019}, where the authors first considered the previous group as a tool for studying latin cycle-sets (not necessarily finite).\\
The displacement group turns out to be also an essential tool to characterize indecomposable q-cycle sets having finite primitive level. In that regard, we exactly compute the primitive level of those having abelian permutation group. We say that a finite indecomposable q-cycle set $X$ has \emph{primitive
level $k$} if $k$ is the biggest positive integer such that there exist q-cycle sets $X_1 = X, X_2, \ldots  ,X_k$ and an epimorphism $p_{i+1}: X_{i}\rightarrow X_{i+1}$  with $|X_i| > |X_{i+1}| > 1$, for every $1\leq  i\leq k-1$, and $X_k$ is primitive. We specify that this notion was originally introduced in \cite{CeOk20x} in the involutive case and it establishes how far a solution is far from being a primitive solution. As a main application, referring to \cite[Question 3.2]{CeOk20x}, we characterize all the involutive solutions of primitive level equal to $2$ and, among these, we completely classify the ones having abelian permutation group.

Finally, we show that the class of multipermutational cycle sets is contained in that of the cycle sets having finite primitive level. Consequently, we pose particular attention to some results and questions arisen in \cite{CeOk20x,rump2020one,smock}.

\bigskip

\section{General results}

This section is devoted to introducing some definitions and results involving the permutation group of a solution, already considered in \cite{rump2019covering} and \cite{CaCaSt20}.
Moreover, we give some notions on the algebraic structure of q-cycle set, that will be useful throughout the paper.
Finally, we extend some results on the retractability of solutions given for the involutive case in \cite{capiru2020} and \cite{RaVe21}.

\medskip

\begin{defin}\label{def:perm-gr-sol}
Let $\left(X,r\right)$ be a solution and consider the permutation $\eta_x$ on $X$ given by $\eta_x\left(y\right):=\rho_{\lambda^{-1}_y\left(x\right)}\left(y\right)$, for all $x,y\in X$. Then, we name the group
\begin{align*}
\mathcal{G}\left(X,r\right):=\seq{\lambda_x, \, \eta_x \ | \,x\in X}
\end{align*}
the \emph{permutation group associated to $\left(X,r\right)$}.
\end{defin}
\noindent Note that if $(X,r)$ is involutive, then $\mathcal{G}(X,r)=\seq{\lambda_x \ | \ x \in X}$, cf. \cite{etingof1998set}.

\medskip

The following theorem characterizes indecomposable solutions in the finite case. We highlight that it is implicitly contained in the paper by Etingof, Soloviev, and Guralnick \cite[Lemma 2.1]{EtSoGu01};
however, to make our exposition self-contained, we give a proof of this result by using our terminology. 

\begin{theor}\label{teoind21}
Let $(X,r)$ be a solution. Then $(X,r)$ is indecomposable if and only if $\mathcal{G}(X,r) $ acts transitively on $X$.
\begin{proof}
To get the claim, we show that $x,y \in X$ are in the same orbit with respect to the action of the group $\mathcal{G}(X,r)$ if and only if they are in the same orbit with respect to the action of the group
\begin{align*}
    \mathcal{F}\left(X,r\right):=\seq{\lambda_x, \, \rho_x \ | \,x\in X}.
\end{align*}
Indeed, in \cite[Proposition 6.6]{CeSmVe19}, it is proved that $\left(X,r\right)$ is indecomposable if and only if $\mathcal{F}\left(X,r\right)$ acts transitively on $X$.\\
Now, suppose that $x$ and $y$ are in the same orbit with respect to the action of $\mathcal{F}\left(X,r\right)$. Then, there exist $n\in \mathbb{N}$ and $x_1,\ldots, x_n\in X$, such that $g_{x_1}\cdots g_{x_n}(x) = y$,
where $g_{x_i}\in \{\lambda_{x_i},\rho_{x_i}\}$, for every $1 \leq i\leq n$. 
We show that there exist $z_1,\ldots, z_n\in X$ and $h_{z_1},\ldots,h_{z_n}\in \mathcal{G}(X,r)$, with $h_{z_i}\in \{\lambda_{z_i},\eta_{z_i}\}$, for every $1 \leq i\leq n$,
such that $h_{z_1}\cdots h_{z_n}(x) = y$, by proceeding by induction on $n$. \\
If $n = 1$ we distinguish two cases. If $g_{x_1}=\lambda_{x_1}$, then we set $z_1:=x_1$ and $h_{z_1}:=\lambda_{z_1}$, while if $g_{x_1}=\rho_{x_1}$, we set $z_1:=\lambda_x\left(x_1\right)$ and $h_{z_1}:=\eta_{z_1}$. \\
Now, we assume that the claim holds for $n-1$. 
Then, there exist $z_2,\ldots,z_{n}\in X$,  $h_{z_2},\ldots,h_{z_{n}}\in \mathcal{G}(X,r)$ such that $h_{z_i}\in \{\lambda_{z_i},\eta_{z_i}\}$, for every $i\in \{2,\ldots,n\}$, and
$$
h_{z_2}\cdots h_{z_n}(x) = g_{x_2}\cdots g_{x_n}(x).
$$
If $g_{x_1}=\lambda_{x_1}$, we set $z_1:= x_1$ and $h_{z_1}:= \lambda_{z_1}$. In this way, we trivially obtain that
$
h_{z_1}\cdots h_{z_n}(x)
=
g_{x_1}\cdots g_{x_n}(x) = y. 
$
Finally, if $g_{x_1}=\rho_{x_1}$, then we set $z_1:=\lambda_{g_{x_2}\cdots g_{x_n}\left(x\right)}\left(x_1\right)$ and $h_{z_1}:=\eta_{z_1}$. Hence,
$$h_{z_1}h_{z_2}\cdots h_{z_n}\left(x\right)=h_{z_1}g_{x_2}\cdots g_{x_n}\left(x\right)=\rho_{\lambda^{-1}_{g_{x_{2}}\cdots g_{x_n}\left(x\right)}\left(z_1\right)}\ g_{x_2}\cdots g_{x_n}\left(x\right)=y,$$
by completing the inductive step. In a similar way, one can show the converse implication.
\end{proof}
\end{theor}

\begin{rem}
In general, even if the orbits of the actions given by $\mathcal{G}(X,r) $ and $\mathcal{F}(X,r) $ coincide, these two groups are different. For instance, if $X:= \{1,2,3,4\}$ and we  consider the following permutations in $\Sym_X$
$$
\lambda_{1}:=(2 \ 3)\qquad \lambda_{2}:=(1\ 4)\qquad \lambda_{3}:=(1\ 2\ 4\ 3)\qquad \lambda_{4}:=(1 \ 3 \ 4 \ 2), 
$$
and we put $\rho_y\left(x\right):=\lambda^{-1}_{\lambda_x\left(y\right)}\left(x\right)$, for all $x,y\in X$, then $(X,r)$ is an involutive solution, cf. \cite[Example 8.2.14]{jespers2007noetherian}. In this case, the group $\mathcal{G}(X,r)$ coincides with $\seq{\lambda_x \ | \ x \in X}$ and it is different from $\mathcal{F}(X,r)$. Indeed, for example, $\rho_1=(2 \ 4)$ does not belong to $\mathcal{G}(X,r)$.
\end{rem}

\medskip

As a consequence of \cref{teoind21}, we obtain that the indecomposability of a solution $(X,r)$ is linked to that of its derived solution, i.e., the pair  $\left(X, r'\right)$ where $r'$ is the map defined by
$
r'\left(x,y\right) = \left(y, \lambda_y\rho_{\lambda^{-1}_x\left(y\right)}\left(x\right)\right)$,
for all $x,y\in X$ (see \cite[p. 579]{So00}).
Clearly, $\mathcal{G}\left(X,r'\right)$ is a subgroup of $\mathcal{G}(X,r)$.
\begin{cor}
Let $\left(X,r\right)$ be a solution and suppose that its derived solution $\left(X,r'\right)$ is indecomposable. Then, $\left(X,r\right)$ is indecomposable.
\begin{proof}
By \cref{teoind21}, $\mathcal{G}(X,r')$ acts transitively on $X$. Since \, $\mathcal{G}(X,r') \leq \mathcal{G}(X,r) $, it follows that $\mathcal{G}(X,r)$ also acts transitively on $X$. Thus, by \cref{teoind21}, the claim follows.
\end{proof}
\end{cor}

\medskip

Now, in light of the existing correspondence between solutions and regular q-cycle sets given in \cite[Proposition 1]{rump2019covering}, let us introduce the analogous notion of indecomposable q-cycle set. To this end, according to \cite[p. 143]{rump2019covering}, we denote by
\begin{align*}
\mathcal{G}\left(X\right)=\seq{\sigma_x, \delta_x \, |\, x\in X}
\end{align*}
the \emph{permutation group associated to a regular q-cycle set $X$}. It is straightforward to check that the group $\mathcal{G}\left(X,r\right)$ coincides with the group $\mathcal{G}\left(X\right)$.

\begin{defin}
A regular q-cycle set $X$ is said to be \emph{indecomposable} if $\mathcal{G}(X)$ acts transitively on $X$. 
\end{defin}
\noindent In view of \cref{teoind21}, a solution $(X,r)$ is indecomposable if and only if the associated regular q-cycle set $X$ is indecomposable.

\medskip

Our aim is to show the retractability of any finite indecomposable q-cycle set $X$ with regular permutation group $\mathcal{G}(X)$. In addition, if $\mathcal{G}(X)$ is abelian, $X$ is multipermutational. To this purpose, we recall the relation of retraction of a regular q-cycle set $X$, contained in \cite[Definition 1]{CaCaSt20}. Specifically, the relation $\sim$  on $X$ given by
$$x\sim y\Longleftrightarrow \sigma_x=\sigma_y 
\quad \mbox{and}\quad   
\delta_x=\delta_y,$$
for all $x,y\in X$, is called the \emph{retract relation} of $X$. Such a relation is a congruence of q-cycle sets and ensures that the quotient $\Ret\left(X\right):=X/\sim$ can be endowed with the two operations induced by $X$, which makes it into a q-cycle set structure. We call such a q-cycle set the \emph{retraction} of $X$.
Analogously, given a solution $(X,r)$, one can define the retract relation of $(X,r)$, by using the maps $\lambda_x$ and $\rho_x$ (see \cite[p. 3595]{JeP18}). As one can expect, the retraction of a solution corresponds to the retraction of a non-degenerate q-cycle set.

\begin{defin}
 A regular q-cycle set $X$ is said to be \emph{retractable} if $|X|=1$ or if there exist two distinct elements $x, y$ of
$X$ such that $\sigma_x=\sigma_y$ and $\delta_x=\delta_y$, otherwise $X$ is \emph{irretractable}. Moreover, $X$ is \emph{multipermutational of level $n$} if $n$ is the minimal non-negative integer such that $|\Ret^n\left(X\right)|=1$, where 
    \begin{align*}\Ret^0\left(X\right):=X \text{ \ and \ } \Ret^i\left(X\right):=\Ret\left(\Ret^{i-1}\left(X\right)\right) \text{, \ for }i>0.\end{align*}
    In this case, we write $\mpl\left(X\right)=n$.
\end{defin}
\noindent Clearly, a multipermutational q-cycle set of level $n$ is retractable, but the converse is not necessarily true. 

\medskip

\begin{lemma}\label{rem_q_q'}
Let $X$ be an indecomposable q-cycle set such that $\mathcal{G}(X)$ is regular. Then, $\mathfrak{q}$ and $\mathfrak{q}'$ coincide if and only if  $\cdot$ and $:$ coincide.
\begin{proof}
Indeed, if $\mathfrak{q}$ and $\mathfrak{q}'$ coincide, then $\sigma_x\left(x\right) = \delta_x\left(x\right)$, for every $x\in X$. Moreover, since $\mathcal{G}(X)$ acts regularly, $\sigma_x = \delta_x$, i.e., the operations $\cdot$ and $:$ coincide. 
The converse implication is trivial.
\end{proof}
\end{lemma}

\medskip

Let us observe that the following result is consistent with those contained in \cite[Proposition 1]{capiru2020} and \cite[Corollary 3.3]{RaVe21} in the context of cycle sets.

\begin{theor}\label{prop_multi}
If $X$ is an indecomposable q-cycle set such that $|X|>1$ and $\mathcal{G}(X)$ is regular, then $X$ is retractable. Moreover, if $\mathcal{G}(X)$ is abelian, then $X$ is multipermutational.
\begin{proof}
If $\mathfrak{q}$ and $\mathfrak{q}'$ coincide, by \cref{rem_q_q'}, $X$ is a cycle set and it is retractable by \cite[Corollary 3.3]{RaVe21}. If $\mathfrak{q}$ and $\mathfrak{q}'$ are different, there exists  $x \in X$ such that $\mathfrak{q}(x)\neq \mathfrak{q}'(x)$. Moreover, by $\eqref{ug1}$, we have that $\sigma_{\mathfrak{q}\left(x\right)}\mathfrak{q}\left(x\right)=\sigma_{\mathfrak{q}'\left(x\right)}\mathfrak{q}\left(x\right)$,
hence, since $\mathcal{G}(X)$ acts regularly on $X$, it follows that $\sigma_{\mathfrak{q}\left(x\right)}=\sigma_{\mathfrak{q}'\left(x\right)}$. With similar arguments, by \eqref{ug2}, we obtain that $\delta_{\mathfrak{q}\left(x\right)}=\delta_{\mathfrak{q}'\left(x\right)}$. Thus, $X$ is retractable.\\
Now, we suppose that $\mathcal{G}(X)$ is abelian and we show the remaining part of the claim by induction on $|X|$. If $|X|=2$, the claim trivially holds. Thus, we assume the claim for every indecomposable q-cycle set having order less than $|X|$ and abelian permutation group. 
Clearly, by the first part of the proof $X$ is retractable and $|\Ret\left(X\right)|<|X|$. 
Furthermore, $\Ret\left(X\right)$ is an indecomposable q-cycle set and  $\mathcal{G}\left(\Ret\left(X\right)\right)$ is again abelian. Therefore, by the inductive hypothesis, $\Ret(X)$ is multipermutational, which completes our claim.
\end{proof}
\end{theor}

\medskip

In the last part of this section, we introduce the following subgroups of the permutation group of a q-cycle set. We underline that they will be essential tools for studying some classes of indecomposable q-cycle sets throughout our paper.
\begin{defin}
Let $X$ be a regular q-cycle set not necessarily finite. Then, the groups
\begin{align*}
\Dis^+(X):=\seq{\sigma_x\sigma_y^{-1},\delta_x\delta_y^{-1}\ |\ x,y\in X}
\qquad
\Dis^-(X):&=\seq{\sigma_x^{-1}\sigma_y,\delta_x^{-1}\delta_y\ |\ x,y\in X}
\end{align*}
are called the \emph{positive displacement group} and the \emph{negative displacement group} of $X$, respectively. 
\end{defin}

\noindent Let us observe that if $X$ is a regular q-cycle set (not necessarily finite), then it holds
$\Dis^{+}(X)\leq \Dis^{-}(X)$. Indeed, by \eqref{ug1} and \eqref{ug2} we have that 
$$
\sigma_{x}\sigma_{y}^{-1}=\sigma_{x\cdot y}^{-1}\sigma_{y:x} \quad \text{and}\quad
\delta_x\delta_y^{-1}=\delta_{x:y}^{-1}\delta_{y\cdot x},
$$
for all $x,y\in X$.
In particular, if $X$ is a q-cycle set in which $\cdot$ and $:$ coincide, then $\Dis^+(X)$ and $\Dis^-(X)$ are exactly the displacement groups introduced by Bonatto, Kinyon, Stanovsk\'{y}, and Vojt\v{e}chovsk\'{y} in \cite[Section 2]{bon2019} in the context of cycle sets.\\
The following results are consistent with Proposition 2.13 and Lemma 2.15 in \cite{bon2019} in the finite case; we underline that, by \cite[Theorem 5]{CaCaSt20}, any finite regular q-cycle set is non-degenerate.

\begin{lemma}\label{bijective}
If $X$ is a regular q-cycle set, the map $\Delta_{(\cdot, \ :)}:X\times X\rightarrow X\times X$ given by 
$$
\Delta_{(\cdot, \ :)}(x,y) = (x\cdot y,y:x),
$$ 
for every $(x,y)\in X\times X$, is bijective.
\begin{proof}
Since $X$ is finite, we only show that $\Delta_{(\cdot, \ :)}$ is injective. 
If $(x,y),(x',y')\in X\times X$ are such that $\Delta_{(\cdot, \ :)}(x,y)=\Delta_{(\cdot, \ :)}(x',y')$, then 
\begin{align*}
    x\cdot y = x'\cdot y' &\Longrightarrow \mathfrak{q}(x\cdot y)=\mathfrak{q}(x'\cdot y')\\
    &\Longrightarrow (y:x)\cdot (y\cdot y)=(y':x')\cdot (y'\cdot y') &\mbox{by \eqref{ug1}}\\
    &\Longrightarrow (y:x)\cdot (y\cdot y)=(y:x)\cdot (y'\cdot y') &\mbox{\text{since} $y:x=y':x'$}\\
    &\Longrightarrow \mathfrak{q}(y) = \mathfrak{q}(y')\\
    &\Longrightarrow y= y'.
\end{align*}
Hence, $y: x =y': x'= y: x'$ and so $x=x'$. Therefore, the claim follows.
\end{proof}
\end{lemma}

\medskip

\begin{prop}\label{dis_dis+_dis-}
Let $X$ be a regular q-cycle set. Then, $\Dis^+(X) = \Dis^-(X)$.
\begin{proof} 
Initially, we know that $\Dis^+(X)\leq \Dis^-(X)$. Moreover, by \cref{bijective}, for every $(x,y)\in X\times X$ there exists a unique pair $(a,b)\in X\times X$ such that $(x,y)=(a\cdot b, b:a)$. Hence, by \eqref{ug1} and \eqref{ug2} we obtain that $\sigma_x^{-1}\sigma^{}_y = \sigma^{}_a\sigma_b^{-1}$ and $\delta_x^{-1}\delta^{}_y = \delta^{}_a\delta_b^{-1}$. Therefore, the claim follows.
\end{proof}
\end{prop}

\noindent In light of \cref{dis_dis+_dis-}, for any regular q-cycle set $X$ we set
$$
\Dis(X):=\Dis^+(X) = \Dis^-(X)
$$
and we call it the \emph{displacement group} of $X$.

\bigskip

\section{Dynamical extensions of indecomposable q-cycle sets and applications}

In this section, we describe indecomposable q-cycle sets in terms of dynamical extensions, and we focus on some applications. At first, we show that the permutation group associated to any indecomposable retractable q-cycle set $X$ always acts  imprimitively on $X$. Moreover, motivated by \cite[Question 9.6(4)]{Ga21}, we give a structure theorem of indecomposable square-free q-cycle sets. Finally, referring to \cite[Question 9.6(7)]{Ga21}, we provide a family of this kind of q-cycle sets that are not self-distributive.

\medskip

Initially, we recall the notion of dynamical pair contained in \cite{CaCaSt20}, that is a useful tool to construct new q-cycle sets. Specifically, given a a q-cycle set $X$, a set $S$, two maps $\alpha:X\times X\times S\to\Sym_S$ and $\alpha':X\times X\times S\to S^S$, where $S^S$ is the set of all the maps from $S$ into itself, the pair $\left(\alpha,\alpha'\right)$ is called a \emph{dynamical pair} if the following equalities
\begin{align}\label{coci}
\alpha_{\left(x\cdot y\right),\left(x\cdot z\right)}\left(\alpha_{\left(x, y\right)}\left(s,t\right),\alpha_{\left(x,z\right)}\left(s,u\right)\right)&=\alpha_{\left(y:x\right),\left(y\cdot z\right)}\left(\alpha'_{\left(y,x\right)}\left(t,s\right),\alpha_{\left( y,z\right)}\left(t,u\right)\right) \notag\\
\alpha'_{\left(x: y\right),\left(x: z\right)}\left(\alpha'_{\left(x,y\right)}\left(s,t\right),\alpha'_{\left(x,z\right)}\left(s,u\right)\right)&=\alpha'_{\left(y\cdot x\right),\left(y:z\right)}\left(\alpha_{\left(y,x\right)}\left(t,s\right),\alpha'_{\left(y,z\right)}\left(t,u\right)\right)\\
\alpha'_{\left(x\cdot y\right),\left(x\cdot z\right)}\left(\alpha_{\left(x,y\right)}\left(s,t\right),\alpha_{\left(x,z\right)}\left(s,u\right)\right)&=\alpha_{\left(y: x\right),\left(y: z\right)}\left(\alpha'_{\left(y,x\right)}\left(t,s\right),\alpha'_{\left(y,z\right)}\left(t,u\right)\right)\notag
\end{align}
hold, for all $x,y,z\in X$ and $s,t,u\in S$.\\ As shown in \cite[Theorem 16]{CaCaSt20}, the triple $(X\times S,\cdot,:)$ where
\begin{align*}
(x,s)\cdot (y,t)
:=(x\cdot y,\alpha_{(x,y)}(s,t)) 
\qquad 
\left(x,s\right): \left(y,t\right):=\left(x: y,\alpha'_{\left(x,y\right)}\left(s,t\right)\right),
\end{align*}
for all $x,y\in X$ and $s,t\in S$, is a q-cycle set. 
If $X$ is regular and $\alpha'_{\left(x,y\right)}\left(s,-\right)\in\Sym_S$, for all $x,y\in X$ and $s\in S$, then  $\left(X\times S,\cdot,:\right)$ is regular. Moreover, the converse is true if $X$ and $S$ have finite order. The q-cycle set $X\times_{\alpha,\alpha'} S:=(X\times S,\cdot,:)$ is said to be a \emph{dynamical extension} of $X$ by $S$.\\
Any q-cycle set can be regarded as a particular dynamical extension by using covering maps. Into detail, given two q-cycle sets $X$ and $Y$, a homomorphism $p:X\rightarrow Y$ is called \textit{covering map} if it is surjective and all the fibers 
$p^{-1}\left(y\right):=\{x\, |\,x \in X,\hspace{1mm} p\left(x\right)=y \}$ have the same cardinality. 
For the ease of the reader, we recall such a result below.
\begin{theor}[Theorem 18, \cite{CaCaSt20}]\label{teo_covering}
Let $X$ and $Y$ be q-cycle sets and $p:X\rightarrow Y$  a covering map. Then, there exist a set $S$ and a dynamical pair $(\alpha,\alpha')$ such that $X\cong Y\times_{\alpha,\alpha'} S$.
\end{theor}

\medskip

The next two results are consistent with those given for cycle sets in  \cite[Section 3]{cacsp2018}.
\begin{lemma}\label{prip}
Every epimorphism $p:X\to Y$ from an indecomposable q-cycle set $X$ to a q-cycle set $Y$ is a covering map.
\begin{proof}
To prove our claim, we show that
\begin{align}\label{pi}
    \sigma_i\left(p^{-1}\left(p\left(x\right)\right)\right)
= p^{-1}\left(p\sigma_i \left(x\right)\right)
\qquad 
 \delta_i \left(p^{-1}\left(p\left(x\right)\right)\right) 
 = p^{-1}\left(p\delta_i \left(x\right)\right),
\end{align}
for all $i, x\in X$. Note that, if $j\in p^{-1}(p(x))$, then we get
$$
p\left(i\cdot j\right)
= p\left(i\right)\cdot p\left(j\right)
= p\left(i\right)\cdot p\left(x\right)
= p\left(i\cdot x\right).
$$ 
Hence, $\sigma_i \left(p^{-1}\left(p\left(x\right)\right)\right)\subseteq p^{-1}\left(p \sigma_i \left(x\right)\right)$.
Conversely, if $j\in p^{-1}\left(p\sigma_i \left(x\right)\right)$, then there exists $k\in X$ such that $j=\sigma_i\left(k\right)$. Thus,
$$
p(i)\cdot p(x) = p(i\cdot x)
= p(j)
= p(i\cdot k)
= p(i)\cdot p(k),
$$
and so $p(x)=p(k)$. Then, they follow that $k\in p^{-1}\left(p\left(x\right)\right)$ and $j\in \sigma_i\left(p^{-1}\left(p\left(x\right)\right)\right)$. In a similar way, one can check the second equality of \eqref{pi}.\\
Now, if $x_1,x_2\in X$, by the indecomposability of $X$, there exists $\pi\in\mathcal{G}\left(X\right)$ such that $\pi(x_1)=x_2$. By the equalities in \eqref{pi}, since $X$ is finite we get $\pi\left(p^{-1}\left(p(x_1\right)\right)=p^{-1}\left(p\left(x_2\right)\right)$. Therefore, we obtain that $| p^{-1}\left(p\left(x_1\right)\right)| = | p^{-1}\left(p\left(x_2\right)\right)|$ and so the map $p$ is a covering map.
\end{proof} 
\end{lemma}

\medskip

\begin{theor}\label{cov}
If $X$ is an indecomposable q-cycle set and $p:X \to Y$ an epimorphism from $X$ to a q-cycle set $Y$, then there exist a set $S$ and a dynamical
pair $(\alpha,\alpha')$ such that $X\cong Y\times_{\alpha,\alpha'} S$. 
\begin{proof}
Initially, by \cref{prip}, the map $p$ is a covering map. Therefore, by \cref{teo_covering}, there exist a set $S$ and a dynamical pair $(\alpha,\alpha')$ such that $X$ is isomorphic to $Y\times_{\alpha,\alpha'} S$, which is our claim.
\end{proof}
\end{theor}

\medskip

As a first application of \cref{cov}, we extend the results obtained in \cite[Remark 2.3 and Lemma 2.4]{CeJeOk20x}, where the authors show that the permutation group $\mathcal{G}(X)$ of a finite retractable cycle set $X$ acts imprimitively on $X$, whenever $X$ has not prime size. To this end, we recall some classical notions of group theory (see, for instance, \cite{dixon1996permutation, dobsonimprimitive, wielandt2014finite}). Given a finite transitive group $G$ acting on a set $X$, then $G$ is said to be \emph{imprimitive} if there exists a subset $\Delta$ of $X$, $\Delta\neq X$, with at least two elements, called \emph{block} for $G$, such that, for each permutation $g$ of $G$, either $g\left(\Delta\right) = \Delta$ or $g\left(\Delta\right)\cap \Delta = \emptyset$. In this way, the set $\{g\left(\Delta\right)\}_{g\in G}$ forms a partition of $X$ which is said to be an \emph{imprimitive blocks system}.

\begin{theor}\label{blocimpr}
Let $X$ be a retractable indecomposable q-cycle set such that $|X|$ is not a prime number. Then, $\mathcal{G}(X)$ acts imprimitively on $X$.
\begin{proof}
If $\mpl(X)=1$, then $\mathcal{G}(X)$ is an abelian group that acts transitively on $X$. Since in this case the action of $\mathcal{G}(X)$ on $X$ is equivalent to the action of $\mathcal{G}(X)$ on itself by left multiplication (see, for example, 1.6.7 in \cite{robinson2012course}), the statement follows from the fact that $\vert X \vert$ is not a prime number.\\
If $\mpl(X)>1$, then $\vert \Ret\left(X\right)\vert >1$ and, by \cref{prip}, the canonical epimorphism $\sigma:X\rightarrow \Ret\left(X\right)$ is a covering map. In this way, by  \cref{cov}, there exist a set $S$ and a dynamical pair $(\alpha,\alpha')$ such that $X$ can be identified with $\Ret\left(X\right)\times_{\alpha,\alpha'} S$. Therefore, we obtain that the partition $\{\{\sigma_x\}\times S \}_{\sigma_x\in \Ret\left(X\right)}$ is an imprimitive blocks system for $X$, hence the claim follows.
\end{proof}
\end{theor}

\medskip

Below, we provide a characterization of indecomposable dynamical extensions of q-cycle sets that includes \cite[Theorem 7]{cacsp2018} given in the context of cycle sets. To this purpose, given a finite q-cycle set $X$, a set $S$, and a dynamical pair $\left(\alpha,\alpha'\right)$, we denote by
\begin{align*}
    H_x:= \Big\{ h \ | \ h\in \mathcal{G}\left(\{x\}\times_{\alpha} S\right),\ h\left(\{x\}\times S\right)= \{x\}\times S \Big\}
\end{align*}
the stabilizer of the set $\{x\}\times S$, for every $x \in X$.

\begin{theor}\label{prop_stabilizer}
Let $X$ be a q-cycle set, $S$ a set, and $\left(\alpha,\alpha'\right)$ a dynamical pair. Then, $X\times_{\alpha,\alpha'} S$ is indecomposable if and only if $X$ is indecomposable and there exists $x\in X$ such that the stabilizer $H_x$ acts transitively on $\{x\}\times S$.
\begin{proof}
Clearly, if $X\times_{\alpha,\alpha'} S$ is indecomposable, then $X$ is indecomposable. Moreover, if $x\in X$ and $s_1,s_2\in S$, by the indecomposability of $X\times_{\alpha,\alpha'} S$, there exists $g\in \mathcal{G}(X\times_{\alpha,\alpha'} S)$ such that $g(x,s_1)=(x,s_2)$. Thus, $g\in H_x$ and the transitivity of $H_x$ on $\{x\}\times S$ follows.\\
Conversely, suppose that $X$ is indecomposable and there exists $x\in X$ such that $H_x$ acts transitively on $\{x\}\times S$. 
By a standard calculation, one can prove that the subgroups $\{H_y\}_{y\in X}$ are all conjugate one to each other. 
Besides, the transitivity of $H_x$ on $\{x\}\times S$ implies the transitivity of $H_y$ on $\{y\}\times S$, for every $y\in X$. 
Since $X$ is indecomposable, if $(x_1,s_1),(x_2,s_2)\in X\times S$, there exist $y_1,\ldots ,y_n\in X$ such that 
$\varphi_{y_1}\cdots\varphi_{y_n}(x_1)=x_2$, where $\varphi_{y_i}\in\{\sigma_{y_i}, \delta_{y_i}\}$, for every $i\in\{1,\ldots, n\}$.
Furthermore, there exist $t_1,\ldots,t_n,\bar{s}$ in $S$ such that $\varphi_{(y_1,t_1)}\cdots\varphi_{(y_n,t_n)}(x_1,s_1)=(x_2,\bar{s})$, where $\varphi_{(y_i,t_i)}\in\{\sigma_{(y_i,t_i)}, \delta_{(y_i,t_i)}\}$, for every $i\in\{1,\ldots, n\}$. Since $H_{x_2}$ acts transitively on $\{x_2\}\times S$, there exists $w\in H_{x_2}$ such that $w(x_2,\bar{s})=(x_2,s_2)$. Therefore, 
$w\varphi_{(x_1,t_1)}\cdots\varphi_{(x_n,t_n)}(x_1,s_1)=(x_2,s_2)$, hence the claim follows.
\end{proof}
\end{theor}

\medskip

\begin{rem}\label{rem_stabilizer}
Note that, assuming that $\mathcal{G}(X\times S)$ acts on $X\times S$, there exists an induced action of $\mathcal{G}(X\times S)$  on $X$ given by $\sigma_{\left(x,s\right)}\left(y\right) = \sigma_x\left(y\right)$  and $\delta_{\left(x,s\right)}\left(y\right) = \delta_x\left(y\right)$, for all $x,y\in X$ and $s\in S$.
Moreover, as a consequence of \cref{prop_stabilizer}, if $X\times_{\alpha, \alpha'}S$ is indecomposable, then
the permutation group $\mathcal{G}(X\times S)$ imprimitively acts on $X \times S$ and the set $\{ \{x\}\times S\}_{x\in X}$ is a blocks system, which we call the \textit{blocks system induced by $X$}.
\end{rem}

\medskip

In the next, we construct examples of indecomposable q-cycle sets by using \cref{prop_stabilizer}.

\begin{ex}\hspace{1mm}
\begin{enumerate}
    \item Let $X:=\mathbb{Z}/4\mathbb{Z}$ be the q-cycle set given by $x\cdot y:=y+1$ and $x:y:=y-1$, for all $x,y\in X$, $S:=\mathbb{Z}/2\mathbb{Z}$, and $\alpha,\alpha':X\times X\times S\mapsto \Sym_S$ the maps given by
    \begin{align*}
        \alpha_{(x,y)}(s,t)=\alpha'_{(x,y)}(s,t)=\begin{cases} t & \mbox{if }x\in \{0,2\} \\ t+1& \mbox{if }x\in \{1,3\}\mbox{ }
\end{cases}
\end{align*}
for all $y\in X$ and $s,t\in S$. Note that $X$ is indecomposable and $\seq{\delta_{(0,0)}\sigma_{(1,0)}}$ is a subgroup of $H_0$ that acts transitively on $\{0\}\times S$. Then, $X\times_{\alpha,\alpha'} S$ is an indecomposable q-cycle set.
\vspace{2mm}
\item Let $k \in \mathbb{N}$, $X:=\mathbb{Z}/2k\mathbb{Z}$ the q-cycle set given by $x\cdot y=x:y:=y+1$, for all $x,y\in X$, $S:=\mathbb{Z}/2\mathbb{Z}$, and $\alpha,\alpha':X\times X\times S\mapsto \Sym_S$ the maps 
    \begin{align*}
        \alpha_{(x,y)}(s,t)=\begin{cases} t & \mbox{if }x\mbox{ is even} \\ t+1& \mbox{if }x \mbox{ is odd}\mbox{ }
\end{cases}
\qquad
        \alpha'_{(x,y)}(s,t)=\begin{cases} t+1 & \mbox{if }x\mbox{ is even} \\ t& \mbox{if }x \mbox{ is odd}\mbox{ }
\end{cases}
\end{align*}
for all $y\in X$ and $s,t\in S$.  Observe that $X$ is indecomposable and $\seq{\sigma_{(0,0)}^{2k-1}\sigma_{(1,0)}}$ is a subgroup of $H_0$ that acts transitively on $\{0\}\times S$. Then, $X\times_{\alpha,\alpha'} S$ is an indecomposable q-cycle set.
\vspace{2mm}
\item Let $p$ be a prime number, $X:=\mathbb{Z}/p\mathbb{Z}$ the q-cycle set given by $x\cdot y=x:y:=y+1$ for all $x,y\in X$, $S:=\mathbb{Z}/p\mathbb{Z}$, and $\alpha,\alpha':X\times X\times S\mapsto \Sym_S$ the maps given by
   \begin{align*}
        \alpha_{(x,y)}(s,t)=t+x\qquad
        \alpha'_{(x,y)}(s,t)=t+x+1,
\end{align*}
for all $x,y\in X$ and $s,t\in S$. Observe that $X$ is indecomposable and $\seq{\sigma_{(0,0)}^{p-1}\sigma_{(1,0)}}$ is a subgroup of $H_0$ that acts transitively on $\{0\}\times S$. Then, $X\times_{\alpha,\alpha'} S$ is an indecomposable q-cycle set.
\end{enumerate}
\end{ex}

\medskip

Inspired by some issues raised by Gateva-Ivanova in \cite{Ga21}, we conclude this section focusing on indecomposable square-free q-cycle sets.
In particular, the author asks for searching specific finite solutions that are indecomposable and square-free.
\begin{que}[Question 9.6(7), \cite{Ga21}]\label{que-GaIv}
Find examples of indecomposable finite square-free solutions which are not left self-distributive.
\end{que}
\noindent We recall that a solution $(X,r)$ is \emph{left self-distributive} if $\rho_x=\id_X$, for every $x\in X$; let us note that such solutions correspond to q-cycle sets for which $\delta_x = \id_X$, for every $x\in X$. For brevity, we call such q-cycle sets left self-distributive. Analogously, one can define the right version of the self-distributivity. For a self-distributive solution or a q-cycle set, we mean a solution or a q-cycle set that is left or right self-distributive.\\
Solutions of the type in \cref{que-GaIv}, that are neither right nor left self-distributive, are unknown to date in the literature and finding instances is challenging. 
Below, we provide a family of such solutions, in terms of dynamical extensions of q-cycle sets, that are indecomposable by \cref{prop_stabilizer}.

\begin{ex}\label{example_1}
Let $S$ be a $2$-elementary abelian group, $X:=\{1,2,3\}$ the indecomposable q-cycle set given by $\sigma_1:=(2\;3)$, $\sigma_2:=(1\;3)$, $\sigma_3:=(1\;2)$ and $\delta_x:=\id_X$, for every $x\in X$, and $\alpha,\alpha':X\times X\times S\rightarrow \Sym_S$ the maps defined by
\begin{align*}
    \alpha_{(x,y)}(s,t)=t
    \qquad\text{and}\qquad
    \alpha'_{(x,y)}(s,t)=\begin{cases} t & \mbox{if }x=y \\ s +t& \mbox{if }x\neq y\mbox{ }
\end{cases},
\end{align*}
for all $x,y\in X$ and $x,t\in S$. 
Then, $X\times_{\alpha,\alpha'} S$ is an indecomposable square-free q-cycle set. 
To show this, we prove the equalities in \eqref{coci}. Initially, note that the first equality is trivially satisfied. 
Moreover, if $x,y,z\in X$ and $s,t,u\in S$, then the third equality reduces to
\begin{align*}
\alpha'_{\left(x\cdot y\right),\left(x\cdot z\right)}\left(t,u\right)=\alpha'_{\left(y,z\right)}\left(t,u\right)   
\end{align*}
and since $y=z$ if and only if $x\cdot y = x\cdot z$, it is satisfied. Now, the second equality reduces to
\begin{align}\label{eq:sec_e}
    \alpha'_{\left(y,z\right)}\left(\alpha'_{\left(x,y\right)}\left(s,t\right),\alpha'_{\left(x,z\right)}\left(s,u\right)\right)=\alpha'_{\left(y\cdot x,z\right)}\left(s,\alpha'_{\left(y,z\right)}\left(t,u\right)\right)
\end{align}
and, to prove it, we need to examine the following cases.
\begin{itemize}
    \item[-] If $x=y=z$, both members in \eqref{eq:sec_e} are equal to $u$.
    \item[-] If $x=y$ and $x\neq z$, then since $X$ is square-free, $y\cdot x=x$ and so $y\cdot x\neq z$. Therefore, the first side in \eqref{eq:sec_e} is equal to $\alpha'_{\left(x,y\right)}\left(s,t\right)+\alpha'_{\left(x,z\right)}\left(s,u\right)=t+\alpha'_{\left(x,z\right)}\left(s,u\right)=t+s+u$ and the second one is equal to $s+\alpha'_{\left(y,z\right)}\left(t,u\right)=s+t+u$.
  \item[-]  If $x\neq y$ and $x=z$, then $y\cdot x\neq x$. Hence, the first member in \eqref{eq:sec_e} is equal to
  $\alpha'_{\left(x,y\right)}\left(s,t\right) + \alpha'_{\left(x,z\right)}\left(s,u\right) = s + t + u$ and the second one is $s+\alpha'_{\left(y,z\right)}\left(t,u\right)=s + t + u$.
  \item[-] If $y = z$ and $x\neq y$, it follows that  $y\cdot x\neq z$, otherwise $z \cdot x=y \cdot x=z=z \cdot z$ and so $x=z$, a contradiction. Thus, the first  side in \eqref{eq:sec_e} is equal to $\alpha'_{\left(x,z\right)}\left(s,u\right)=s+u$ and the second one is equal to $s+\alpha'_{\left(y,z\right)}\left(t,u\right)=s+u$.
  \item[-] If $x, y, z$ are all distinct, by the definition of the maps $\sigma_t$, we have that $y\cdot x = z$. Hence, the first member in \eqref{eq:sec_e} is  $\alpha'_{\left(x,y\right)}\left(s,t\right)+\alpha'_{\left(x,z\right)}\left(s,u\right)=s+t+s+u=u+t$ and the second one is equal to $\alpha'_{\left(y,z\right)}\left(t,u\right)=t+u$.
\end{itemize}
By \cite[Theorem 16]{CaCaSt20}, $X\times_{\alpha,\alpha'} S$ is a q-cycle set. To show that it is indecomposable, by \cref{prop_stabilizer}, it is sufficient to prove that $H_0$ acts transitively on $\{0\}\times S$. This fact is true, since we have that  $\delta_{(1,s)}(0,0)=(1:0,\,\alpha'_{(1,0)}(s,0))=(0,s)$, for every $s\in S$. Moreover, it is easy to check that $X\times_{\alpha,\alpha'} S$ is square-free.
Finally, let us observe that, if $s\in S$, with $s\neq 0$, then
\begin{align*}
    \sigma_{(1,s)}(2,0) = (3,0)
    \qquad\text{and}\qquad
    \delta_{(1,s)}(2,0) = (2, s),
\end{align*}
hence $\sigma_{(1,s)}\neq\id_{X\times S}$ and $\delta_{(1,s)}\neq\id_{X\times S}$.
Therefore, $X\times_{\alpha,\alpha'} S$ is neither right nor left self-distributive.
\end{ex}

\medskip

Let us highlight that the smallest q-cycle set obtained by using the technique in \cref{example_1} is of order $6$.
By computer calculations, we checked that there are not indecomposable square-free q-cycle sets of order $2$, $3$, and $4$, that are not self-distributive. At present, we are not able to state if there exist q-cycle sets of order $5$ of such a type.

\medskip

Finally, we give a structure theorem for indecomposable square-free q-cycle sets in terms of dynamical extensions which partially goes in the direction of \cite[Question 9.6(4)]{Ga21}.

\begin{lemma}\label{prop_squarefree}
Let $X$ be an indecomposable square-free q-cycle set. Then, $\Ret^k(X)$ is a q-cycle set in which $\cdot$ and $:$ do not coincide, for every $k\in \mathbb{N}$. In particular, $X$ is not a multipermutational q-cycle set.
\begin{proof}
At first, observe that clearly $\Ret^k\left(X\right)$ is indecomposable, for every $k\in \mathbb{N}$. Now, suppose that there exists $\bar{k}\in \mathbb{N}$ such that $\Ret^{\bar{k}}\left(X\right)$ is a q-cycle set in which the operations $\cdot$ and $:$ are equal. Then, by \cite[Theorem 1]{rump2005decomposition}, $\Ret^{\bar{k}}\left(X\right)$ is decomposable, a contradiction.\\
Moreover,  if $X$ is multipermutational of level $h$, since $X$ is square-free, then $\Ret^{h-1}\left(X\right)$ is the trivial decomposable cycle set given by $x\cdot y = y$, for all $x,y\in \Ret^{h-1}\left(X\right)$, an absurd.
\end{proof}
\end{lemma}

\medskip

\begin{theor}
Let $X$ be an indecomposable square-free q-cycle set. Then, there exist a set $S$, an irretractable square-free q-cycle set $Y$, and a dynamical pair $\left(\alpha,\alpha'\right)$ such that  $X\cong Y\times_{\alpha,\alpha'} S$.
\begin{proof}
By \cref{prop_squarefree}, there exists $k\in \mathbb{N}$ such that $\Ret^k\left(X\right)$ is an irretractable square-free q-cycle set. Therefore, by \cref{cov}, there exist a set $S$ and a dynamical pair $\left(\alpha,\alpha'\right)$ such that  $X$ is isomorphic to $\Ret^k\left(X\right) \times_{\alpha,\alpha'} S$.
\end{proof}
\end{theor}

\bigskip

\section{Simple regular q-cycle sets}

This section aims to give a group-theoretic characterization of simple regular q-cycle sets, of which we illustrate some applications. Moreover, we show that any simple regular q-cycle set is indecomposable.

\medskip

The notion of simplicity for q-cycle sets has been already  introduced in \cite[Definition 3]{CaCaSt20} and is consistent with that given by
Vendramin for cycle sets in \cite[Definition 2.9]{vendramin2016extensions}. However, in the context of cycle sets, Ced\'o and Okni\'nski gave another definition in \cite[Definition 3.3]{CeOk20x}, which coincides with that provided by Vendramin in the indecomposable case. 
For our purposes, we consider the next definition which includes that in \cite{CeOk20x}.
\begin{defin}\label{def-simple}
   Let $X$ be a regular q-cycle set with $|X| > 1$. Then, $X$ is called \emph{simple} if, for every epimorphism $f : X \to Y$ of q-cycle sets, either $f$ is an isomorphism or $|Y | = 1$.
\end{defin}

\medskip

By \cref{prip}, if $X$ is an indecomposable q-cycle set, \cref{def-simple} coincides with \cite[Definition 3]{CaCaSt20}. Let us note the fact that $X$ is indecomposable is not surprising since this is the case of every simple q-cycle set, exactly as it happens for cycle sets, cf. \cite[Lemma 4.1]{CeOk20x}.

\begin{prop}\label{prop_simple_inde}
Let $X$ be a simple regular q-cycle set with $|X|>1$. Then, $X$ is indecomposable.
\begin{proof}
    Assume that $X$ is decomposable. It follows that there exists a non-trivial orbit $A$ with $|A|<|X|$. 
    Let us consider the trivial q-cycle set on the set $Y:= \{1,2\}$, i.e., $x\cdot y = x:y = y$, for all $x,y\in Y$, and the surjective map $f:X\to Y$ given by
    \begin{align*}
        f\left(x\right)
        := \begin{cases}
        1 &\text{if} \ x\in A\\
        2 &\text{if} \ x\in X\setminus A
        \end{cases}
    \end{align*}
    and we show that $f$ is a homomorphism of q-cycle sets. Let $x \in X$.
    If $y\in A$, since $A$ is an orbit we have that $x\cdot y\in A$, thus $f\left(x\cdot y\right) = 1 = f\left(y\right) = f\left(x\right)\cdot f\left(y\right)$. 
    Now, if $y\in X\setminus A$, clearly $x\cdot y\in X\setminus A$, thus  $f\left(x\cdot y\right) = 2 = f\left(y\right) = f\left(x\right)\cdot f\left(y\right)$. 
    Similarly, $f\left(x\right):f\left(y\right) = f\left(x:y\right)$, for all $x,y\in X$. Therefore, we obtain a contradiction and so the claim follows.
\end{proof}
\end{prop}

\medskip 

Now, we give the characterization of all simple regular q-cycle sets which we mentioned before. Using the terminology contained in \cite{dixon1996permutation, dobsonimprimitive}, given an imprimitive blocks system $\mathcal{B}:=\{g\left(\Delta\right)\}_{g \in G}$ for $G$, if $\Delta\in\mathcal{B}$, the \emph{set-wise stabilizer} of $\Delta$ in $G$ is the subgroup of $G$ given by
\begin{align*}
G_{\Delta}:=\{g\ |\ g\in G,\  g(\Delta)=\Delta  \}.
\end{align*}
\noindent Moreover, the \emph{fixer} of $\mathcal{B}$ is the subgroup of $G$ given by  
$$
\Fix(\mathcal{B})=\bigcap\limits_{\Delta \in \mathcal{B}} G_{\Delta}. 
$$
If $X$ is a q-cycle set, it is useful recalling that if $\mathcal{G}(X)$ is imprimitive and $\{\Delta_x\}_{x\in X}$ is an imprimitive blocks system, then  $g\left(\Delta_x\right) = \Delta_{g(x)}$, for all $g\in \mathcal{G}(X)$ and $x\in X$.
\medskip

\begin{lemma}\label{lemma_simple}
Let $X$ be a regular q-cycle set. If $\mathcal{G}\left(X\right)$ is imprimitive on $X$ and $\{\Delta_x\}_{x\in X}$ is an imprimitive blocks system, the relation defined by
\begin{align}\label{doppia_tilde}
    \forall \ x,y\in X\qquad x\approx y\Longleftrightarrow \Delta_x = \Delta_y
\end{align}
is a left congruence of $X$. 
\begin{proof}
Clearly, $\approx$ is an equivalence relation on $X$. Moreover, if $y,z\in X$ satisfy $y\approx z$, we get
\begin{align*}
 y\approx z
 &\Longrightarrow \sigma_x(\Delta_y)=\sigma_x(\Delta_z)
 \quad \text{and}\quad \delta_x(\Delta_y) = \delta_x(\Delta_z)\\
 &\Longrightarrow\Delta_{x\cdot y} = \Delta_{x\cdot z}
 \quad \text{and}\quad\Delta_{x : y} = \Delta_{x : z}\\
 &\Longrightarrow x\cdot y\approx x\cdot z  \qquad \text{and} \qquad x : y\approx x : z,
\end{align*}
for every $x\in X$.
\end{proof}
\end{lemma}

\medskip

In the following, given a regular q-cycle set $X$, if $\Delta$ is a block for $\mathcal{G}(X)$, we denote by $\Dis(X, \Delta)$ the subgroup of $\Dis(X)$ given by
\begin{align*}
   \Dis(X,\Delta):=\seq{\sigma_x^{-1}\sigma_y,\delta_x^{-1}\delta_y \ |\ x,y\in \Delta}.
\end{align*}

\medskip

\begin{theor}\label{th:char-simple}
Let $X$ be a regular q-cycle set such that $|X| > 1$. Then, $X$ is simple if and only if $\mathcal{G}(X)$ acts transitively on $X$ and, for every non-trivial imprimitive blocks system $\{\Delta_x\}_{x\in X}$, there exists 
$u\in X$ such that $\Dis(X, \Delta_u)\nleq   \Fix\left(\{\Delta_x\}_{x\in X} \right)$.
\begin{proof}
Suppose that $\mathcal{G}(X)$ acts transitively on $X$ and, for every non-trivial imprimitive blocks system $\{\Delta_x\}_{x\in X}$, there exists $\Delta\in \{\Delta_x\}_{x\in X}$ such that $\Dis(X, \Delta)\nleq \Fix(\{\Delta_x\}_{x\in X})$. If $X$ is not simple, then there exist a q-cycle set $Y$ and a non trivial covering map $p:X\rightarrow Y$. Thus, if $x\in X$, set
$$
\Delta_x:=\{z\in X \ |\ p(x)=p(z)\},
$$
we prove that $\{\Delta_x\}_{x\in X}$ is an imprimitive blocks system for $\mathcal{G}(X)$. To this end, without loss of generality, we consider $y \in X$ and $\sigma_y \in \mathcal{G}(X)$ such that $\sigma_y(\Delta_x) \cap \Delta_x \neq \emptyset$ and it is enough to prove that $\sigma_y(\Delta_x)\subseteq \Delta_x$. Thus, if $z\in \Delta_x$, since there exists $t\in\Delta_x$ such that $t = \sigma_y\left(a\right)$,  for a certain $a\in \Delta_x$, we have that
\begin{align*}
    p(\sigma_y\left(z\right))
    = p\left(y \cdot z\right)
    = p\left(y \cdot a\right) 
    = p\left(t\right) 
    = p\left(x\right),
\end{align*}
i.e., $\sigma_y\left(z\right) \in \Delta_x$. If $u\in X$, it follows that $\sigma_x\left(\Delta_u\right) = \Delta_{x\cdot u}$.
Moreover, if $x'\in X$ is such that $\Delta_x = \Delta_{x'}$, since
$$
p(x\cdot u) = p(x)\cdot p(u) = p(x')\cdot p(u) = p(x'\cdot u),
$$
we obtain that $\sigma_x\left(\Delta_u\right) = \sigma_{x'}\left(\Delta_u\right)$.
Similarly, one can show that $\delta_x(\Delta_u) = \delta_{x'}(\Delta_u)$. Clearly, these facts imply that $\Dis\left(X, \Delta_x\right)\leq \mathcal{G}(X)_{\Delta_u}$, for every $u\in X$, that is an absurd.\\
Conversely, assume that $X$ is simple. Then, by \cref{prop_simple_inde}, $X$ is indecomposable. By contradiction, suppose that $\mathcal{G}(X)$ has a non-trivial imprimitive blocks system $\{\Delta_x\}_{x\in X}$ such that $\Dis(X,\Delta_u)\leq \Fix\left(\{\Delta_x\}_{x\in X}\right)$, for every $u\in X$. To get the claim, since $X$ is finite, we show that the relation $\approx$ in \eqref{doppia_tilde}
is a congruence. In this way, by \cref{lemma_simple}, the canonical map $p:X\rightarrow X/\approx$ is an epimorphism of finite q-cycle sets and it is not trivial since $\{\Delta_x\}_{x \in X}$ is a non-trivial imprimitive blocks system, a contradiction. Thus, let $x,x'\in X$  such that $x\approx x'$. Then, $\Delta_x=\Delta_{x'}$ and, if $y\in X$, since  $\Dis(X,{\Delta_x})\leq \mathcal{G}(X)_{\Delta_y}$ 
it follows that $\sigma_x^{-1}\sigma^{}_{x'}(\Delta_y)=\Delta_y$ and $\delta_x^{-1}\delta^{}_{x'}(\Delta_y) = \Delta_y$. Hence, $\sigma_{x'}(\Delta_y)=\sigma_x(\Delta_y)$ and  $\delta_{x'}\left(\Delta_y\right) = \delta_x\left(\Delta_y\right)$ and so $ \Delta_{x\cdot y} = \Delta_{x'\cdot y}$ and $\Delta_{x : y} = \Delta_{x':\,y}$. Therefore, $x\cdot y\approx x'\cdot y$ and $\quad x:y \approx x':y$, which is the desired conclusion.
\end{proof}
\end{theor}

\medskip
\begin{rem}\label{quoziente_doppiatilde}
If $X$ is an indecomposable q-cycle set with $|X| > 1$ and $\mathcal{G}(X)$ has a non-trivial imprimitive blocks system $\{\Delta_x\}_{x\in X}$ such that $\Dis(X,\Delta_u)\leq \Fix\left(\{\Delta_x\}_{x\in X}\right)$, for every $u\in X$, then $\{\Delta_x\}_{x\in X}$ can be equipped with a q-cycle set structure. Indeed, by the proof of \cref{th:char-simple}, one can observe that the quotient $X/\approx$ coincides with $\{\Delta_x\}_{x\in X}$.
\end{rem}

\medskip

Now, we close this section by showing how to use concretely \cref{th:char-simple} on specific q-cycle sets. 

\begin{ex} (cf.  \cite[Example 2.11]{vendramin2016extensions})
Let $X$ be the indecomposable q-cycle set on $\{1,2,3,4\}$ given by \begin{center}
 $\sigma_1=\delta_1:=(1 \ 4)$, \, $\sigma_2=\delta_2:=(1 \ 3 \ 4 \ 2)$, \, $\sigma_3=\delta_3:=(2 \ 3)$, \, $\sigma_4=\delta_4:=(1 \ 2 \ 4\ 3)$.   
\end{center} Then, by \cref{th:char-simple}, $X$ is simple. Indeed, the only imprimitive blocks system is given by $\{\Delta_1,\Delta_2\}$, where $\Delta_1:=\{1,4\}$ and $\Delta_2:=\{2,3\}$, and, since $\sigma_1\sigma_4^{-1}(\Delta_1)=\Delta_2$, clearly $\Dis(X,\Delta_1)\nleq \mathcal{G}(X)_{\Delta_1}$.
\end{ex}

\medskip

\begin{ex} (cf. \cite[Remark 4.11]{CeOk20x})
Let $X$ be the indecomposable q-cycle set on $\{1,2,3,4,5,6,7,8,9\}$ given by
\begin{align*}
   &\sigma_1=\delta_1:=(1\ 3\ 8\ 4\ 5\ 2\ 9\ 7\ 6)
   \qquad 
   &\sigma_2=\delta_2:= (1\ 7\ 6\ 4\ 3\ 8\ 9\ 5\ 2)\\
   &\sigma_3=\delta_3:= (1\ 7\ 8\ 4\ 3\ 2\ 9\ 5\ 6)
   \qquad
   &\sigma_4=\delta_4:=(1\ 2\ 7\ 4\ 6\ 3\ 9\ 8\ 5)\\
   &\sigma_5=\delta_5:= (1\ 8\ 5\ 4\ 2\ 7\ 9\ 6\ 3)
   \qquad
   &\sigma_6=\delta_6:= (1\ 8\ 7\ 4\ 2\ 3\ 9\ 6\ 5)\\
   &\sigma_7=\delta_7:=(1\ 9\ 4)(2\ 8\ 6)
   &\sigma_8=\delta_8:= (1\ 9\ 4)(3\ 7\ 5)\qquad\\
   &\sigma_9 =\delta_9:= (2\ 8\ 6)(3\ 7\ 5).
\end{align*}
Then, by \cref{th:char-simple}, $X$ is simple. In fact, the only imprimitive blocks system is $\{\Delta_1, \, \Delta_2, \, \Delta_3\}$, where $\Delta_1:=\{1,4,9\}$, $\Delta_2:=\{2,6,8\}$, and $\Delta_3:=\{3,5,7 \}$. Note that $\Dis(X,\Delta_1)\nleq \mathcal{G}(X)_{\Delta_1}$ since $\sigma_1\sigma_9^{-1}(\Delta_1)=\Delta_2$.
\end{ex}

\medskip

\begin{ex}
Let $X:=\{1,2,3,4,5,6 \}$ be the indecomposable q-cycle set given by
\begin{align*}
   \sigma_1:=(2\ 4 \ 5 \ 3)
   \qquad\sigma_2:= (1 \ 3 \ 6 \ 4)\\
   \qquad\sigma_3:= (1\ 5 \ 6 \ 2) 
   \qquad\sigma_4:=(1\ 2\ 6 \ 5) \\
   \qquad\sigma_5:= (1\ 4 \ 6 \ 3)
   \qquad\sigma_6:= (2 \ 3 \ 5 \ 4)
\end{align*}
and $\delta_x:=\id_X$, for every $x\in X$. The only imprimitive blocks system is $\{\Delta_1, \, \Delta_2, \, \Delta_3\}$, where $\Delta_1:=\{1,6\}$, $\Delta_2:=\{2,5\}$, and $\Delta_3:=\{3,4 \}$. Note that $\Dis(X,\Delta_u)=\mathcal{G}(X)_{\Delta_v}$, for all $u,v\in X$. Hence, by \cref{th:char-simple}, $X$ is not simple.
\end{ex}

\bigskip

\section{Primitive level of indecomposable q-cycle sets}

In this section, we introduce the notion of primitive level of q-cycle sets, which is consistent with the one given in \cite[p. 7]{CeOk20x} for cycle sets. Moreover, we characterize all the indecomposable q-cycle sets having a finite primitive level in group-theoretic terms. In addition, we compute the primitive level of those having abelian permutation group.

\medskip

\begin{defin}
A q-cycle set $X$ is said to be \emph{primitive} if $\mathcal{G}\left(X\right)$ acts primitively on $X$. Moreover, we say that a finite indecomposable q-cycle set $X$ has \emph{primitive
level $k$} if $k$ is the biggest positive integer such that 
\begin{enumerate}
    \item[(1)] there exist q-cycle sets $X_1 = X, X_2, \ldots  ,X_k$, with $|X_i| > |X_{i+1}| > 1$, for every $1\leq  i\leq k-1$;
    \item[(2)] there exists an epimorphism of q-cycle sets $p_{i+1}: X_{i}\rightarrow X_{i+1}$, for every $1\leq  i\leq k-1$;
    \item[(3)] $X_k$ is primitive.
\end{enumerate}
\end{defin}

\noindent Observe that, by \cref{prip}, if $X$ is an indecomposable q-cycle set of order $p_1^{\alpha_1}\cdots p_n^{\alpha_n}$, for some primes $p_1,\ldots,p_n$, and with primitive level $k$, then $k$ is at most $\alpha_1+\cdots+\alpha_n$. 
Clearly, q-cycle sets having primitive level $1$ are exactly those for which $\mathcal{G}(X)$ acts primitively on $X$. Cycle sets with primitive level $1$ are completely classified (cf. \cite[Theorem 3.1]{CeJeOk20x} and \cite[Theorem 2.13]{etingof1998set}): they are, up to isomorphisms, the ones on $\mathbb{Z}/p\mathbb{Z}$ (where $p$ is a prime number) given by $x\cdot y = y + 1$, for all $x,y\in \mathbb{Z}/p\mathbb{Z}$.\\
In general, there exist primitive q-cycle sets whose cardinality is not a prime number, as one can see in the following example.

\begin{ex}
Let $X:=\{1,2,3,4\}$ be the q-cycle set given by \begin{align*}
 \sigma_1:=(2\;4\;3) \qquad \sigma_2:=(1\;3\;4)\qquad \sigma_3:=(1\;4\;2)\qquad 
 \sigma_4:=(1\;2\;3),    
\end{align*}
and $\delta_x:=\id_X$, for every $x\in X$. Then, $X$ is a primitive q-cycle set.
\end{ex}

\medskip

On the other hand, not all q-cycle sets have finite primitive level: for example, simple q-cycle sets with imprimitive permutation group have not finite primitive level. In the next, we give a characterization of all q-cycle sets having finite primitive level, in which the displacement group and its subgroups have a crucial role. 
We recall that if $G$ is a group that acts transitively on a finite set $X$, then a block $\Delta$ of $X$ is said to be \emph{maximal} if, for every block $\Delta'$ such that $\Delta\subseteq \Delta'\subseteq X$, it follows that $\Delta=\Delta'$ or $\Delta'=X$ (cf. \cite{wielandt2014finite}). 
Moreover, an imprimitive blocks system $\mathcal{B}$ is said to be maximal if one (and hence all) of its blocks is maximal. Clearly, an imprimitive blocks system $\mathcal{B}$ is maximal if and only if the induced action of $G$ on $\mathcal{B}$ is primitive.

\begin{theor}\label{liv_finito}
Let $X$ be an indecomposable q-cycle set. Then, $X$ has finite primitive level if and only if there exists a maximal imprimitive blocks system $\{\Delta_x\}_{x\in X}$ such that $\Dis(X,\Delta)\leq \Fix(\{\Delta_x\}_{x\in X})$, for every $\Delta \in \{\Delta_x \}_{x\in X}$. 
\begin{proof}
If $X$ has finite primitive level $k$, there exist q-cycle sets $X_1=X, X_2, \dots  ,X_k$, an epimorphism $p_{i+1}: X_{i}\rightarrow X_{i+1}$, with $|X_i| > |X_{i+1}| > 1$, for every $1\leq  i\leq k-1$, and $X_k$ is primitive. Clearly, the composition $p_kp_{k-1}\cdots p_2$ is an epimorphism from $X$ to $X_k$, hence, by \cref{teo_covering}, there exist a set $S$ and a dynamical pair $(\alpha,\alpha')$ such that $X$ can be identified with $X_k\times_{\alpha,\alpha'} S$. Since $X_k$ is primitive, the action of $\mathcal{G}\left(X_k\times_{\alpha,\alpha'} S\right)$ on $X_k$ (as in \cref{rem_stabilizer}) is primitive. It follows that the imprimitive blocks system $\{\{x\}\times S\}_{x\in X_k}$ induced by $X_k$ is maximal (see \cite[p. 18]{wielandt2014finite}). 
Moreover, if $(x_1,s_1)$ and $(x_2,s_2)$ are in the same block, it follows that $x_1 = x_2$ and 
$$
\sigma_{\left(x_1,s_1\right)}^{-1}\sigma^{}_{\left(x_2,s_2\right)}\left(\{y\}\times S\right)
= \{y\}\times S, 
$$
for every $y\in X_k$. Analogously, one can see that $\delta_{\left(x_1,s_1\right)}^{-1}\delta^{}_{\left(x_2,s_2\right)}\left(\{y\}\times S\right)
= \{y\}\times S$. Hence, the necessary condition follows.\\
To prove the vice versa, it is sufficient to find an epimorphism from $X$ to a primitive q-cycle set. Thus, if $\{\Delta_x\}_{x\in X}$ is a maximal imprimitive blocks system such that $\Dis(X,\Delta)\leq \Fix(\{\Delta_x\}_{x\in X})$, for every $\Delta \in \{\Delta_x \}_{x\in X}$, then, by  \cref{th:char-simple},
$X$ is not simple. Furthermore, by \cref{quoziente_doppiatilde}, $\{\Delta_x\}_{x\in X}$ has a structure of q-cycle set such that $p:X\rightarrow \{\Delta_x\}_{x\in X}$, $x\mapsto \Delta_x$ is an epimorphism of q-cycle sets. Hence, by \cref{teo_covering}, there exist a set $S$ and a dynamical pair $(\alpha,\alpha')$ such that $X$ can be identified with $\{\Delta_x\}_{x\in X} \times_{\alpha,\alpha'} S$. Moreover, since $\{\Delta_x\}_{x\in X}$ is a maximal imprimitive blocks system, we have that the action of $\mathcal{G}(X)$ on $\{\Delta_x\}_{x\in X}$ is primitive, and this clearly implies that $\mathcal{G}(\{\Delta_x\}_{x\in X}) $ acts primitively on $\{\Delta_x\}_{x\in X}$. Therefore, $p$ is the requested epimorphism. 
\end{proof}
\end{theor}

\medskip

Now, our aim is to show how to compute the finite primitive level of any indecomposable q-cycle set $X$ having $\mathcal{G}(X)$ as an abelian group.
Let us recall that a transitive action of a finite abelian group $G$ on a finite set $X$ is equivalent to the left regular action of $G$ (for more details see, for instance, 1.6.7 in \cite{robinson2012course}). For this reason, without loss of generality, if $X$ is an indecomposable q-cycle set having abelian permutation group, we can assume that $X$ is an abelian group and, for every $x\in X$, there exist $h_x, h'_x\in X$ such that $\sigma_x$ and $\delta_x$ are the translations $t_{h_x}$ and $t_{h'_x}$ by $h_x$ and $h_{x}'$, respectively, namely
\begin{align*}
    \sigma_x\left(y\right) =  h_x + y
    \qquad
    \delta_x\left(y\right) = h'_x + y,
\end{align*}
for every $y \in X$. Hereinafter, if $H$ is a subgroup of $G$, we denote by $\sim_H$ the usual equivalence relation with respect to $H$, i.e., $x\sim_H y$ if and only if $x-y\in H$, for all $x,y\in G$.

\medskip

\begin{lemma}\label{fin1}
    Let $X$ be an indecomposable q-cycle and assume that $\mathcal{G}(X)$ is abelian. 
   Then, there exists a non-trivial subgroup $H$ of $X$ such that $\Ret(X) = X/H$. In addition, if $H'$ a subgroup of $X$ contained in $H$,  the equivalence relation induced by $H'$ is a congruence of q-cycle sets.
    \begin{proof} Note that, by \cref{prop_multi}, $X$ is retractable. If $|\Ret(X)|=1$, put $H:= X$, the first part of our claim follows. Otherwise, since $X$ has not prime order, by \cref{blocimpr}, the retraction induces an imprimitive blocks system $\mathcal{B}$. By a standard calculation, denoted by $\Delta_0 \in \mathcal{B}$ the block containing $0$, one has that $H:=\Delta_0$ is a subgroup of $X$. It follows that $\Ret(X)=\{x+H\}_{x\in X}$.\\
    Now, let $x,x',y,y'\in X$ such that $x\sim_{H'}x'$ and $y\sim_{H'}y'$. Note that, since $H'\subseteq H$, $\sigma_x = \sigma_{x'}$ and $\delta_x = \delta_{x'}$. Hence, there exists $h\in X$ such that $\sigma_x$ and $\sigma_{x'}$ are the translations by $h$. Thus, we get $y+h\sim_{H'} y'+h$ and so, $x\cdot y\sim_{H'} x'\cdot y'$. Similarly, one can show that $x: y\sim_{H'} x': y'$. Therefore, the claim follows.
    \end{proof}
\end{lemma}

\medskip

\begin{theor}\label{calcolo_prim_ab}
    Let $X$ be an indecomposable q-cycle set with $|X|>1$ and such that $\mathcal{G}(X)$ is abelian and has size $p_1^{\alpha_1}\cdots p_m^{\alpha_m}$, with $p_1, \ldots, p_m$ primes. Then, $X$ has primitive level equal to $\alpha_1+\cdots+\alpha_m$.
\begin{proof}
    We proceed by induction on $|X| = |\mathcal{G}(X)|$. If $|X| = 2$ the statement is trivial. Now, suppose that the claim is true for every q-cycle set $Y$ with $|Y|<|X|$ and having abelian permutation group. By \cref{fin1}, there exists a non-trivial subgroup $H$ of $X$ such that $\Ret(X)=X/H$. Without loss of generality, we can suppose that $p_1 \mid |H|$. Thus, let $H'$ be a subgroup of $H$ such that $|H'| = p_1$. Then, by \cref{fin1}, $X/H'$ is a q-cycle set and $p:X\rightarrow X/H'$ is an epimorphism of q-cycle sets. Finally, by the inductive hypothesis on $X/H'$, we have that the primitive level of $X/H'$ is $\alpha_1 + \cdots + \alpha_{m}-1$. Therefore, the primitive level of $X$ is at least  $\alpha_1+\cdots+\alpha_m$ and hence the statement is proved.
\end{proof}
\end{theor}

\medskip

As a consequence of \cref{calcolo_prim_ab}, the class of multipermutational indecomposable q-cycle sets is contained in that of finite primitive level.
\begin{cor}\label{inclus}
    Let $X$ be a multipermutational indecomposable q-cycle set. Then, $X$ has finite primitive level.
\begin{proof}
Assume that $X$ has multipermutational level $n$ and let $p: X \to \Ret^{n-1}(X)$ be the canonical epimorphism. Then, $\Ret^{n-1}(X)$ is a q-cycle set having abelian permutation group. Hence, if $\Ret^{n-1}(X)$ has prime size, the claim is proved. Otherwise, by \cref{calcolo_prim_ab}, there exists an epimorphism $\bar{p}$ from $\Ret^{n-1}(X)$ to a primitive q-cycle set $Y$, with $|Y| < |\Ret^{n-1}(X)|$. Therefore, $\bar{p}p$ is an epimorphism from $X$ to $Y$ and the statement follows. 
\end{proof}
\end{cor}

\bigskip

\section{Some applications to the indecomposable involutive solutions}
 
In this section, we specialize some of the
results previously obtained to indecomposable cycle sets. Into the specific, following \cite[Questions 3.2 and 7.3]{CeOk20x}, we characterize indecomposable cycle set having primitive level $2$, and we give a sufficient condition to guarantee the simplicity of indecomposable cycle sets of size $p^2$. 
 Moreover, we discuss some matters contained in \cite{CeOk20x,rump2020one, smock}.

\medskip

The following is a characterization of indecomposable cycle sets having finite primitive level. Unlike \cref{liv_finito}, one can note that this result entirely involves the displacement group.

 \begin{theor}\label{liv_finito_cyc}
 Let $X$ be an indecomposable cycle set. Then, $X$ has finite primitive level if and only if there exists an imprimitive blocks system $\{\Delta_x\}_{x\in X}$ having prime size such that $\Dis(X)\leq \Fix(\{\Delta_x\}_{x\in X})$.
\begin{proof}
If $X$ has finite primitive level, with a similar argument in the proof of \cref{liv_finito}, we obtain that $X$ can be identified with a dynamical extension $X_k\times_{\alpha,\alpha'} S$, where $X_k$ is a primitive cycle set. By  \cite[Theorem 3.1]{CeJeOk20x}, $|X_k|=p$, with $p$ a prime number, and $\sigma_{x}(y)=\gamma(y)$, for all $x,y\in X_k$, where $\gamma$ is a $p$-cycle. These facts imply that 
\begin{center}
    $\sigma_{(x_1,s_1)}^{-1}\sigma^{}_{(x_2,s_2)}\left(\{x\}\times S\right)=\{x\}\times S$,
\end{center} 
for all $x \in X_k$ and $(x_1,s_1),(x_2,s_2)\in X_k\times S$. Hence, if we set $\Delta_{(x,s)}:=\{x\}\times S$, for every $(x,s)\in X_k\times S$, we obtain that the imprimitive blocks system $\{\Delta_{(x,s)}\}_{(x,s)\in X_k\times S}$ has size $p$ and $\Dis(X)\leq \Fix(\{\Delta_{(x,s)}\}_{(x,s)\in X_k\times S})$.\\
Conversely, since $\Dis(X, \Delta)\leq \Fix(\{\Delta_x\}_{x\in X})$, for every $\Delta \in \{\Delta_x\}_{x\in X}$, and
every imprimitive blocks system $\{\Delta_x\}_{x\in X}$ having prime size is maximal, the vice versa is a consequence of \cref{liv_finito}.
\end{proof}\end{theor}
\medskip

Among involutive solutions of finite primitive level, Ced{\'o} and Okni{\'n}ski asked the following.
\begin{que}(\cite[Question 3.2]{CeOk20x})
Describe involutive solutions of primitive level $2$.
\end{que}

\medskip

Using \cref{liv_finito_cyc}, we obtain a characterization of indecomposable cycle sets having primitive level $2$. To this purpose, we recall that, given a group $G$ which acts transitively on a finite set $X$ and two imprimitive blocks systems $\mathcal{B},\mathcal{B}'$ of $G$, then $\mathcal{B}'$ is a \textit{refinement} of $\mathcal{B}$ if there exist $\Delta\in \mathcal{B}$ and $\Delta'\in \mathcal{B}' $ such that $\Delta' \subseteq \Delta$.

\begin{theor}\label{prop-lev2}
Let $X$ be an indecomposable cycle set such that $|X|$ is not prime. Then, $X$ has primitive level $2$ if and only if the following conditions hold:
\begin{enumerate}
    \item there exists an imprimitive blocks system $\{\Delta_x\}_{x\in X}$ having prime size such that $\Dis\left(X\right)\leq \Fix\left(\{\Delta_x\}_{x\in X}\right)$;
    \item for every non-trivial refinement $\{\Delta'_x\}_{x\in X}$ of an imprimitive blocks system $\{\Delta_x\}_{x\in X}$ satisfying $1.$, there exists $\Delta'\in \{\Delta'_x\}_{x\in X}$ such that $\Dis\left(X,\Delta'\right)\nleq \Fix\left(\{\Delta'_x\}_{x\in X}\right)$.
\end{enumerate}
\begin{proof}
If $X$ has primitive level $2$ then, by \cref{liv_finito_cyc}, the condition $1.$ is satisfied. Thus, let $\{\Delta_x\}_{x\in X}$ be such an imprimitive blocks system having prime size and such that $\Dis\left(X\right)\leq \Fix\left(\{\Delta_x\}_{x\in X}\right)$ and consider a non-trivial refinement  $\{\Delta'_x\}_{x\in X}$ such that $\Dis(X,\Delta')\leq \Fix(\{\Delta'_x\}_{x\in X})$, for every $\Delta'$ in $\{\Delta'_x\}_{x\in X}$. It follows that the maps $p_1,p_2$ given by $p_1:X\rightarrow \{\Delta'_x\}_{x\in X}$, $x\mapsto \Delta'_x$ and $p_2:\{\Delta'_x\}_{x\in X}\rightarrow \{\Delta_x\}_{x\in X}$, $\Delta_x'\mapsto \Delta_x$ are epimorphisms of cycle sets, hence $X$ has primitive level at least $3$, a contradiction.\\
Conversely, if conditions $1.$ and $2.$ hold, then, by \cref{liv_finito_cyc}, $X$ has finite primitive level. Moreover, if $X$ has primitive level greater than $2$, there exist two indecomposable cycle sets $X_1$ and $X_2$, where $X_2$ is primitive, and two epimorphisms $p_1:X\rightarrow X_1$ and $p_2:X_1\rightarrow X_2$.
Hence, by \cref{teo_covering}, $X_1$ can be identified with a cycle set having $X_2\times T$ as underlying set, for a suitable set $T$. In the same way, $X$ can be identified with a cycle set having $X_1\times S=X_2\times T\times S$ as underlying set, for a suitable set $S$. Therefore, the set $\{\{a\}\times S\}_{a\in X_1}$ is a refinement of  $\{\{u\}\times T\times S\}_{u\in X_2}$. Moreover, by  \cref{th:char-simple}, $\Dis\left(X,\{a\}\times S\right)\leq \Fix\left(\{\{x\}\times S\}_{x\in X_1}\right)$, for every $\{a\}\times S\in \{\{x\}\times S\}_{x\in X_1}$, a contradiction.
\end{proof}\end{theor}

\medskip

As a consequence of \cref{calcolo_prim_ab}, we concretely describe cycle sets having primitive level $2$ and abelian permutation group, up to isomorphism.
\begin{cor}
    All the indecomposable cycle sets with abelian permutation group and primitive level $2$ have size $pq$ with $p,q$ not necessarily distinct primes. Hence, they are the ones provided in \cite[Theorem 21]{capiru2020}.
 \begin{proof}
    It follows by \cref{calcolo_prim_ab} and \cite[Theorem 21]{capiru2020}.
 \end{proof}
\end{cor} 

In general, by \cref{prip}, the retractable indecomposable cycle sets having size $pq$, with $p,q$ not necessarily distinct primes, have primitive level equal to $2$.
In \cite{CeOk20x}, one can find the next question:
\begin{que}[Question 7.5, \cite{CeOk20x}]\label{que7.5}
Does there exist a simple involutive solution $(X, r)$ such that $|X|=p_1p_2 \cdots p_n$, for $n >1$, with $p_1, p_2, \dots, p_n$ distinct primes?
\end{que}

\noindent Inspecting all the indecomposable cycle sets having size $6$ (using a GAP package in \cite{Ve15pack}), we note that they are all retractable. These facts clearly imply that it does not exist a simple indecomposable cycle set having size $6$, answering in this specific case in the negative sense to the \cref{que7.5}.

\medskip
 
In \cite[Corollary 6.6]{smock}, the authors give a sufficient condition for the decomposability of a multipermutational finite cycle set. For the ease of the reader, we recall below such a result.

\begin{cor}\cite[Corollary 6.6]{smock}\label{cor_smok}
Let $X$ be a multipermutational cycle set with $|X|>1$ and suppose that, for every $x \in X$, there exists $y\in X$ such that $x\cdot y = y$ and $y\cdot x = x$. Then, $X$ is decomposable.
\end{cor}

In \cite[Question 6.7]{smock} Smoktunowicz and Smoktunowicz asked the following:
 \begin{que}
    Is \cref{cor_smok} also true without the assumption that the multipermutational level of $\left(X, r\right)$ is finite?
 \end{que}
 
In \cite[Theorem 4]{rump2020one}, Rump showed that \cref{cor_smok}, in general, is not true if the cycle set is not of multipermutational type; nevertheless the hypotheses can be relaxed. In this context, we show that \cref{cor_smok} is not unexpected: indeed, using dynamical extensions and \cref{inclus}, such cycle sets belong to a larger class of decomposable cycle sets. 
 
 \begin{prop}\label{epi}
Let $X,Y$ be cycle sets, $p:X\rightarrow Y$ a covering map, and suppose that $Y$ is an indecomposable cycle set having prime size. Then, $x\cdot y\neq y$, for all $x,y\in X$.
\begin{proof}
If $x,y\in X$ are such that $x\cdot y = y$, then $p(x)\cdot p(y) = p(x\cdot y) = p(y)$, which contradicts \cite[Theorem 2.13]{etingof1998set}.
\end{proof}
\end{prop}

\medskip

Below, we provide method, alternative to \cref{cor_smok}, to check that a cycle set is decomposable without the hypothesis of multipermutationality. 
\begin{cor}\label{corepi}
    Let $X,Y$ be cycle sets, with $|Y|$ a prime number. Suppose that $Y$ is an epimorphic image of $X$ and that there exist $x,y\in X$ such that $x\cdot y = y$. Then, $X$ is decomposable.
\begin{proof}
If by absurd $X$ is indecomposable, then so $Y$ is. Hence, by \cref{prip} and \cref{epi}, we obtain a contradiction. Therefore, the statement is proved.
\end{proof}\end{cor}

\medskip

Now, let us  focus on two applications of \cref{corepi}: in the first one, we give a direct simple proof of \cref{cor_smok} which does not make use of the brace theory; in the second one, we extend \cite[Theorem 4]{rump2020one} in the case of finite cycle sets, giving more information on the cyclic structure of the left multiplications.

\begin{proof}[Proof of \cref{cor_smok}]
If $X$ is indecomposable, by \cref{inclus}, it has finite primitive level. Hence, there exists an epimorphism $p$ from $X$ to an indecomposable cycle set $Y$ having prime size. Therefore the claim follows by \cref{corepi}.
\end{proof}

\medskip
The next result clearly implies \cite[Theorem 4]{rump2020one} in the finite case.

\begin{cor}\label{cor:fixp}
    Let $X$ be an indecomposable cycle set having finite primitive level. Then, $\sigma_x$ has no fixed point, for every $x\in X$.
\begin{proof}
Since $X$ has finite primitive level, there exists an epimorphism from $X$ to a cycle set $Y$ having prime size. If there exists $x\in X$ such that $\sigma_x$ has a fixed point, by \cref{corepi}, $X$ is decomposable, a contradiction.
\end{proof}
\end{cor}

\medskip

Let us note that \cref{cor:fixp} allows for giving a sufficient condition for \cite[Question 7.3]{CeOk20x}, which we recall below.
\begin{que}
Let $p$ be a prime and $(X, r)$ an indecomposable and irretractable solution of cardinality $p^2$. Is $(X, r)$ simple? 
\end{que}

\medskip

\begin{cor}\label{cor_ceok_answer}
Let $X$ be an indecomposable cycle set having size $p^2$, for some prime $p$. If there exist $x,y\in X$ such that $x\cdot y = y$, then $X$ is simple.
\begin{proof}
Clearly, $X$ is simple or $X$ has primitive level $2$. In the second case, by  \cref{cor:fixp}, it follows that $x\cdot y\neq y$, for all $x,y\in X$, against the hypothesis. 
\end{proof}
\end{cor}

\medskip

\noindent Finally, observe that the converse of \cref{cor_ceok_answer} does not hold: indeed, in \cite[p. 19]{CeOk20x} the authors give a simple indecomposable cycle set of size $9$, where the left multiplications have no fixed points.

\bigskip


\bibliographystyle{elsart-num-sort}
\bibliography{Bibliography}

\end{document}